\documentclass{amsart}

\usepackage{graphicx, color}

\usepackage{amsthm}
\usepackage{amssymb}
\usepackage{amsmath}
\usepackage{graphicx}
\usepackage{amscd}
\usepackage{amsfonts}
\usepackage{amsbsy}
\usepackage[T1]{fontenc}
\usepackage[english]{babel}

\newtheorem{innercustomgeneric}{\customgenericname}
\providecommand{\customgenericname}{}
\newcommand{\newcustomtheorem}[2]{%
  \newenvironment{#1}[1]
  {%
   \renewcommand\customgenericname{#2}%
   \renewcommand\theinnercustomgeneric{##1}%
   \innercustomgeneric
  }
  {\endinnercustomgeneric}
}

\newcustomtheorem{customthm}{Theorem}
\usepackage{epsfig}
\usepackage{amssymb}

\newtheorem*{cor}{Corollary}

\newtheorem{ex}{Example}

\title[A measure-theoretic expansion exponent]{A measure-theoretic expansion exponent}

\author{C.A. Morales}

\address{Beijing Advanced Innovation Center for Future Blockchain and Privacy Computing, Beihang University, Beijing, China \& Beijing Academy of Blockchain and edge Computing, Beijing, 100086.}
\email{morales@impa.br}

\keywords{Expansion exponent, Borel probability measure, Metric space}

\subjclass[2020]{Primary 37B25; Secondary 37B65}

\begin{document}

\begin{abstract}
The expansion exponent (or expansion constant) for maps was introduced by Schreiber in \cite{s}.
In this paper, we introduce the analogous exponent for measures. We shall prove the following results:
The expansion exponent of a measurable map is equal to the minimum of the expansion exponents taken over all Borel probability measures. In particular,  
a map expands small distances (in the sense of Reddy \cite{r}) if and only if every Borel probability measure has a positive expansion exponent.
Any nonatomic invariant measure with a positive expansion exponent is positively expansive in the sense of \cite{m}.  
For ergodic invariant measures, the Kolmogorov–Sinai entropy is bounded below by the product of the expansion exponent and the measure’s upper capacity.  
As a consequence, any ergodic invariant measure with both positive upper capacity and positive expansion exponent must have positive entropy.
\end{abstract}

\maketitle

\section{Introduction}
\noindent
Let \( T: X \to X \) be a measurable map on a metric space.  
This means that the preimage of any Borel set under \( T \) is also a Borel set—that is, an element of the Borel sigma-algebra generated by the open subsets of \( X \).
The {\em expansion exponent} (or expansion constant) of \( T \), denoted \( E(T) \), is the supremum of all \( \lambda \in \mathbb{R} \) for which there exists \( \epsilon > 0 \) such that  
\(
d(T(x), T(y)) \geq e^{\lambda} d(x, y)
\)
for all \( x, y \in X \) with \( d(x, y) < \epsilon \).
This exponent was studied by Schreiber \cite{s} for smooth maps of compact Riemannian manifolds.

In this paper, we introduce an analogous exponent for Borel probability measures. These are $\sigma$-additive measures $\mu$ defined in the Borel sigma-algebra with $\mu(X)=1$.
For such measures we define the expansion exponent $E_\mu(T)$ with respect to $T$
as the supremum of all $\lambda\in\mathbb{R}$ for which there exists $\epsilon>0$ such that the event $d(T(x),T(y))<e^\lambda d(x,y)$ given $d(x,y)<\epsilon$ is negligible for all $x\in X$. In other words,
\begin{multline*}
E_\mu(T)=\sup\bigg\{\lambda\in\mathbb{R}\mid \exists \epsilon>0\mbox{ such that }
\mu\{\{y\in B(x,\epsilon)\mid\\
 d(T(x),T(y))<e^\lambda d(x,y)\})=0,\forall x\in X\bigg\}
\end{multline*}
(with the convention $\sup\emptyset=-\infty$).

This exponent measures the strongest uniform expansion rate at which nearby points (in the sense of \( \mu \)) are expanded under \( T \), in a measure-theoretic rather than pointwise or differentiable sense. The condition requires that, for each \( x \in X \), the set of \( \mu \)-measure points around \( x \) that fail to expand at rate \( \lambda \) under \( T \) has zero measure. The term “expansion exponent” is motivated by its role as a global, measure-theoretic analogue of the classical Lyapunov expansion exponent, extending the notion of local expansion exponential divergence to non-differentiable or non-smooth settings via probabilistic concentration.

We apply this measure-theoretical exponent by pproving the following results:
The expansion exponent of a measurable map is the minimum of the expansion exponents over the Borel probability measures. In particular, 
a map expands small distances (in the sense of Reddy \cite{r}) if and only if every Borel probability measure has a positive expansion exponent. 
Furthermore, any nonatomic invariant measure with positive expansion exponent is positively expansive in the sense of \cite{m}. 
Finally, for ergodic invariant measures that the product of the expansion exponent and the measure's upper capacity is bounded above by the Kolmogorov–Sinai entropy of the measure. 
As a consequence, every ergodic invariant measure with positive upper capacity and positive expansion exponent has positive entropy.

Before the precise statements we present some few examples.
For the first one recall that the {\em Dirac measure} supported on $p\in X$ is the Borel probability measure $\delta_p$ defined by $\delta_p(B)=1$ or $0$ depending on whether $p\in B$ or not, for all Borelian $B$. A measure is Dirac if it is the Dirac measure supported on some point.

\begin{ex}
\label{point}
Let $T:X\to X$ be a measurable map of a metric space.
Then, $\delta_p$ is has positive expansion exponent if and only if
$p$ is a {\em locally expanding point}, i.e.,
there are $k>1$ and $\epsilon>0$ such that $d(T(p),T(x))\geq kd(p,x)$ for all $x\in B(p,\epsilon).$ From this it is easy to construct continuous maps with Dirac measures with positive expansion exponent.
\end{ex}

We say that a Borel probability measure $\mu$ is the {\em convex combination of Dirac measures} if there are sequences $(x_i)_{i\in\mathbb{N}}$ and
$(t_i)_{i\in\mathbb{N}}$ with $0\leq t_i$ and $\sum_{i\in\mathbb{N}}t_i=1$ such that
$$
\mu=\sum_{i\in\mathbb{N}}t_i\delta_{x_i}.
$$

\begin{ex}
If $T:X\to X$ is constant (i.e. $T(x)=T(y)$ for all $x,y\in X$), then
$E_\mu(T)\neq-\infty$ for some Borel probability measure $\mu$ $\iff$ $X$ is {\em uniformly $\mu$-discrete} i.e. there is $\epsilon>0$ such that $\mu(B(x,\epsilon))=\mu(\{x\})$ for all $x\in X$. In particular, if $X$ is separable, then $\mu$ is the convex combination of Dirac measures.
\end{ex}

More examples can be obtained as follows.

\begin{ex}
The following properties hold for all measurable map $T:X\to X$, all Borel probability measures $\mu,\nu$ of $X$ and $0<t<1$:
\begin{enumerate}
\item
$E_{t\mu+(1-t)\nu}(T)=\min\{E_\mu(T),E_\nu(T)\}$ (in particular, the map $\mu\mapsto E_\mu(T)$ is convex);
\item
If $\mu\prec\nu$ (i.e. $\nu(A)=0$ $\implies$ $\mu(A)=0$),
then $E_\mu(T)\geq E_\nu(T)$.
\end{enumerate}
It follows that the set of Borel probability measures with expansion exponent bigger than certain $\lambda$, i.e., if $\mu$ and $\nu$ have both expansion exponent bigger than $\lambda$, then
so does $t\mu+(1-t)\nu$ for all $0\leq t\leq 1$.
\end{ex}

Recall that a metric space $X$ is {\em Lindel\"of} if every open cover has a countable subcover. A map $T:X\to X$ is termed {\em distance preserving} if $d(T(x),T(y))= d(x,y)$ for all $x,y\in X$. Obviously, every distance-presrving map is continuous hence measurable.
Among those maps we can mention the identity of $X$ and the circle rotations.

\begin{ex}
\label{exo}
If $T:X\to X$ is a distance-preserving map of a Lindel\"of metric space, then $E_\mu(T)\leq0$ for all Borel probability measure $\mu$.
\end{ex}

\begin{proof}
If $\mu$ were a measure with positive expansion exponent of $T$,
then it would exists $k>1$ and $\epsilon>0$ satisfying \eqref{pito}
hence $0=\mu(\{y\in B(x,\epsilon)\mid d(T(x),T(y))<kd(x,y)\})=\mu(B(x,\epsilon))$ for all $x\in X$. Since the open covering $\{B(x,\epsilon):x\in X\}$ has a countable subcover
(for $X$ is Lindel\"of), we would have that $\mu(X)=0$ which is absurd.
\end{proof}

Now, we present our results.
The first one is that the expansion exponent of a map is the infimum of the
measure-theoretical expansion exponent. More precisely, we have the following result (compare with Theorem 1 in \cite{s}).

\begin{customthm}{A}
\label{white}
If $T:X\to X$ is a measurable map of a metric space,
\begin{equation}
\label{chugule}
E(T)=\min_\mu E_\mu(T),
\end{equation}
where the minimum is over the Borel probability measures $\mu$ of $X$.
\end{customthm}

A corollary in the spirit of \cite{lmn} is as follows.
Recall that a map $T:X\to X$ {\em expands small distances} \cite{r} if
there are $k>1$ and $\epsilon>0$ such that
$$
d(T(x),T(y))\geq kd(x,y),\quad\quad\forall x\in X,\, y\in B(x,\epsilon).
$$
Here $B(x,\epsilon)$ stands for the open $\epsilon$-ball with center $x$ while $B[x,\epsilon]$ will denote the corresponding closed ball.

\begin{cor}
A measurable map of a metric space $T:X\to X$ expands small distances if and only if every Borel probability measure has positive expansion exponent.
\end{cor}

For another application, we recall that a map of a metric space $T:X\to X$ is {\em positively expansive} if
there is $\epsilon>0$ such that $x=y$ whenever $x,y\in X$ and $d(T^i(x),T^i(y))\leq \epsilon$ for all integer $i\geq 0$.

It is known that a continuous map of a compact metric space is positively expansive if and only if it expands small distances up to some equivalent metric \cite{r}.
This motivate us to relate measures with positive expansion exponent with the positively expansive measures \cite{m}.

On the other hand, a Borel probability measure $\mu$ is {\em positively expansive} for a measurable map $T:X\to X$ if there is $\epsilon>0$ (called {\em expansivity constant}) such that
$\mu(\Phi_\epsilon(x))=0$ for all $x\in X$ where
$$
\Phi_\epsilon(x)=\{y\in X:d(T^i(x),T^i(y))\leq\epsilon,\,\forall i\geq0\}.
$$

Clearly, every positively expansive measure is {\em nonatomic} i.e. $\mu(\{x\})=0$ for all $x\in X$.
Since there are Borel probability measures with positive expansion exponent which are not nonatomic (Example \ref{point}), there are also Borel probability measures with positive expansion exponent which are not positively expansive.
We now give some sufficient conditions for a nonatomic Borel probability measure with positive expansion exponent to be positively expansive. Recall that a Borel probability measure $\mu$ is {\em invariant} if $\mu(T^{-1}(B))=\mu(B)$ for all Borelian $B\subset X$.

Then, we will prove the following result.

\begin{customthm}{B}
\label{walter}
Every nonatomic invariant measure with positive expansion exponent of a measurable map of a metric space is positively expansive.
\end{customthm}

We say that $x\in X$ is an {\em equicontinuous point} of a map $T:X\to X$ if
for every $\epsilon>0$ there is $\delta>0$ such that if $y\in X$ and $d(x,y)<\delta$, then
$d(T^i(x),T^i(y))<\epsilon$ for all integer $i\geq0$.

As in Example 2.5 of \cite{m} it can be proved that there are no positively expansive measures for pointwise equicontinuous maps of Lindel\"of metric spaces.
Then, we obtain the following corollary from Theorem \ref{walter} related to Example \ref{exo}.

\begin{cor}
A pointwise equicontinuous map of a Lindel\"of metric space does not have nonatomic invariant measures with positive expansion exponent.
\end{cor}

To state our last result we need the following definitions.
The first is the {\em upper capacity} \cite{dzg} of a Borel probability measure $\mu$ of $X$ defined by
\begin{equation}
\label{faraday}
\overline{dim}_B(\mu)=\lim_{\delta\to0}\limsup_{\beta\to0}\frac{\log N_\mu(X,\beta,\delta)}{-\log \beta},
\end{equation}
where $N_\mu(X,\beta,\delta)$ is the minimal number of open $\beta$-balls needed to cover a subset of $\mu$-measure bigger than or equal to $1-\delta$.
This definition is due to Ledrappier \cite{l}.

Next, given a collection $\mathcal{X}_0,\cdots,\mathcal{X}_{n-1}$ of subsets of $X$ we write
$$
\bigvee_{i=0}^{n-1}\mathcal{X}_i=\{A_0\cap\cdots A_{n-1}\mid A_i\in \mathcal{X}_i,\, \forall i=0,\cdots, n-1\}.
$$
We also define
$T^{-i}(\mathcal{X}_0)=\{T^{-i}(A)\mid A\in \mathcal{X}_0\}$ for all integer $i\geq0$ and all $T:X\to X$.

A {\em partition} of $X$ is a collection of Borel subsets $\alpha$ whose reunion is $X$.
Clearly, if $T$ is measurable, then
$T^{-i}(\alpha)$ is a partition for all $i\geq0$.
If $\mu$ is a Borel probability measure of $X$, then we define the entropy of a partition $\alpha$ with respect to $\mu$ as
$$
H_\mu(\alpha)=-\sum_{A\in \alpha}\mu(A)\log\mu(A)
$$
(with the convention $0\cdot\log0=0$).

The (Kolmogorov-Sinai) entropy \cite{w} of an invariant measure $\mu$ (w.r.t $T$) is defined by
$$
h_\mu(T)=\sup_\alpha\lim_{n\to\infty}\frac{1}n H(\bigvee_{i=0}^{n-1}T^{-i}(\alpha))
$$
where the supremum is over the finite partitions $\alpha$.

Finally, we say that a Borel probability measure $\mu$ is {\em ergodic} for a measurable map $T:X\to X$
if $\mu(A)\in \{0,1\}$ for all Borelian $A\subset X$ such that $T^{-1}(A)=A$. See \cite{w} for more details.

With these definitions we have the following result. It can be compared with Theorem 1 in \cite{dzg} for expansive homeomorphisms.

\begin{customthm}{C}
\label{ermantraut}
If $T:X\to X$ is a continuous map of a compact metric space, and $\mu$ is an ergodic invariant measure,
\begin{equation}
\label{jimmy}
h_\mu(T)\geq \overline{dim}_B(\mu)\cdot E_\mu(T).
\end{equation}
\end{customthm}

It is known that for every ergodic invariant measure $\mu$
of a continuous map of a compact metric space $T:X\to X$ if
$h_\mu(T)>0$, then $\mu$ is positively expansive \cite{am}.
The converse fails namely there are ergodic invariant measures $\mu$ for certain continuous maps of compact metric spaces $T:X\to X$ such that $\mu$ is positively expansive but $h_\mu(T)=0$
(e.g. the Denjoy circle maps \cite{m}).
Nevertheless, a kind of converse is given by the following corollary.

\begin{cor}
If $T:X\to X$ is a continuous map of a compact metric space, and $\mu$ is an ergodic invariant measure of $T$ with positive expansion exponent and positive upper capacity, then
$h_\mu(T)>0$.
\end{cor}

\section{Proof of the theorems}

\noindent

\begin{proof}[Proof of Theorem \ref{white}]
Let $T:X\to X$ be a measurable map of a metric space.
First we show
\begin{equation}
\label{chugu}
E(T)=\inf_\mu E_\mu(T),
\end{equation}
where the infimum is over the Borel probability measures $\mu$.

If $\lambda\in\mathbb{R}$ and
there is $\epsilon>0$ such that $d(T(x),T(y))\geq e^\lambda d(x,y)$ for all $x,y\in X$ with
$d(x,y)<\epsilon$, then
$\mu(\{y\in B(x,\epsilon)\mid d(T(x),T(y))<e^\lambda d(x,y)\})=\mu(\emptyset)=0$ for all $x\in X$ and all Borel probability measure $\mu$.
Therefore,
$\lambda<E_\mu(T)$ for all Borel probability measure $\mu$ hence
$$
E(T)\leq \inf_\mu E_\mu(T).
$$

To prove the converse inequality assume by contradiction that it is not true.
Then, we can choose $\lambda$ such that
\begin{equation}
\label{palo}
E(T)<\lambda<\inf_\mu E_\mu(T).
\end{equation}
The first inequality implies that there are sequences $x_i,y_i\in X$ such that
$$
d(x_i,y_i)<\frac{1}i\quad\quad\mbox{ and }\quad\quad d(T(x_i),T(y_i))<e^\lambda d(x_i,y_i),\quad\quad\forall i\in\mathbb{N}.
$$
Define
$$
\mu_*=\sum_{i=1}^\infty \frac{1}{2^{i}}\delta_{y_i}.
$$
Then, $\mu$ is a Borel probability measure of $X$ so $E_{\mu_*}(T)>\lambda$ by the sEond inequality in \eqref{palo}. Thus,
there is $\epsilon>0$ such that
$$
\mu_*(\{y\in B(x,\epsilon)\mid d(T(x),T(y))<e^\lambda d(x,y)\})=0,\quad\quad\forall x\in X.
$$
Now, take $j\in\mathbb{N}$ large so that $\frac{1}j<\epsilon$.
Then,
$$
\mu_*(\{y\in B(x_j,\epsilon)\mid d(T(x_j),T(y))<e^\lambda d(x_j,y)\})=0.
$$
However, $d(x_j,y_j)<\frac{1}j<\epsilon$ and $d(T(x_j),T(y_j))<e^\lambda d(x_j,y_j)$
thus $y_j\in \{y\in B(x_j,\epsilon)\mid d(T(x_j),T(y))<e^\lambda d(x_j,y)\}$
hence
\begin{eqnarray*}
\frac{1}{2^{j}}&\leq& \left(\sum_{i=1}^\infty\frac{1}{2^{i}}\delta_{y_i}\right)(\{y_j\})\\
&=& \mu_*(\{y_j\})\\
&\leq& \mu_*(\{y\in B(x_j,\epsilon)\mid d(T(x_j),T(y))<kd(x_j,y)\})\\
&=&0
\end{eqnarray*}
which is absurd.
This contradiction proves
$$
\inf_\mu E_\mu(T)\leq E(T)
$$
hence \eqref{chugu} holds.

Next, we observe that
$$
E_{\sum_{i\in\mathbb{N}} t_i\mu_i}(T)\leq \inf_{i\in\mathbb{N}}E_{\mu_i}(T),
$$
for all sequence of Borel probability measures $(\mu_i)_{i\in\mathbb{N}}$ and all sequence of positive numbers $(t_i)_{i\in\mathbb{N}}$ with $\sum_{i\in\mathbb{N}}t_i=1$.
In fact, if $\lambda\in\mathbb{R}$ and $\epsilon>0$ satisfy
$$
\left(
\sum_{i\in\mathbb{N}} t_i\mu_i\right)
(\{y\in B(x,\epsilon)\mid d(T(x),T(y))<e^\lambda d(x,y)\})=0,\quad\quad\forall x\in X,
$$
then $\mu_i(\{y\in B(x,\epsilon)\mid d(T(x),T(y))<e^\lambda\})=0$ for all $x\in X$ and $i\in\mathbb{N}$ hence
$\lambda<E_{\mu_i}(T)$ for all $i\in\mathbb{N}$ proving the result.

By \eqref{chugu} we can choose a sequence of probabilities $\mu_i$ such that
$E_{\mu_i}(T)\to E(T)$ as $i\to\infty$.
Define
$$
\mu=\sum_{i\in\mathbb{N}} \frac{1}{2^i}\mu_i
$$
which is a Borel probability measure measure so
$E_\mu(T)\geq E(T)$ by \eqref{chugu}. Moreover, the above observation implies
$$
E_\mu(T)\leq \inf_{i\in\mathbb{N}} E_{\mu_i}(T)= E(T).
$$
Then, $E_\mu(T)=E(T)$ which together with \eqref{chugu} proves
\eqref{chugule}.
This ends the proof.
\end{proof}

In the sequel we shall use the following remark:
A Borel probability measure $\mu$ has positive expansion exponent with respect to a measurable map $T:X\to X$ if and only if there are $k>1$ and $\epsilon>0$ satisfying
\begin{equation}
\label{pito}
\mu(\{y\in B(x,\epsilon)\mid d(T(x),T(y))<kd(x,y)\})=0,\quad\quad\forall x\in X.
\end{equation}

\begin{proof}[Proof of Theorem \ref{walter}]
Let $\mu$ be a nonatomic invariant measure of a measurable map of a metric space $T:X\to X$. Suppose that $\mu$ has positive expansion exponent.
Then, there are $k>1$ and $\epsilon>0$ satisfying \eqref{pito}.
So, defining
\begin{equation}
\label{tos}
C(x)=\{y\in B(x,\epsilon)\mid d(T(x),T(y))\geq kd(x,y)\},\quad\quad\forall x\in X,
\end{equation}
we obtain
\begin{equation}
\label{pa}
\mu(B(x,\epsilon)\setminus C(x))=0,\quad\quad\forall x\in X.
\end{equation}

We assert that
\begin{equation}
\label{pajuo}
\Phi_{\frac{\epsilon}2}(x)\setminus \{x\}\subset \bigcup_{i\geq0}T^{-i}(B(T^i(x),\epsilon)\setminus C(T^i(x))),\quad\quad\forall x\in X.
\end{equation}
Otherwise, there are $x\in X$ and $y\in \Phi_{\frac{\epsilon}2}(x)\setminus \{x\}$ such that
$$
y\notin \bigcup_{i\geq0} T^{-i}(B(T^i(x),\epsilon)\setminus C(T^i(x))).
$$
Since $y\in \Phi_{\frac{\epsilon}2}(x)$, one has $y\in T^{-i}(B(T^i(x),\epsilon))$ for all $i\geq0$ then
$$
y\in \bigcap_{i\geq0}T^{-i}(C(T^i(x))).
$$
Now, for any $n\in \mathbb{N}$ the above implies $T^{n-1}(y)\in C(T^{n-1}(x))$
so
$$
d(T^{n}(x),T^{n}(y))\geq kd(T^{n-1}(x),T^{n-1}(y)).
$$
Likewise, $T^{n-2}(y)\in C(T^{n-2}(x))$ so
$$
d(T^{n-1}(x),T^{n-1}(y))\geq kd(T^{n-2}(x),T^{n-2}(y)).
$$
Repeating this process we obtain
$$
d(T^{n-i}(x),T^{n-i}(y))\geq kd(T^{n-i-1}(x),T^{n-i-1}(y)),\quad\quad\forall i=0,\cdots n-1.
$$
It follows that
\begin{multline*}
d(T^{n}(x),T^{n}(y))\geq kd(T^{n-1}(x),T^{n-1}(y))
\geq k^2d(T^{n-2}(x),T^{n-2}(y))\geq\\
\geq \cdots\geq k^{n}d(x,y).
\end{multline*}
But $y\in \Phi_{\frac{\epsilon}2}(x)$ so $d(T^{n}(x),T^{n}(y))\leq\frac{\epsilon}2$ hence
$$
d(x,y)\leq \frac{\epsilon}{2k^n},\quad\quad\forall n\geq0.
$$
By letting $n\to\infty$ we get $d(x,y)=0$ namely $x=y$ which contradicts that $y\in\Phi_{\frac{\epsilon}2}(x)\setminus \{x\}$.
This proves \eqref{pajuo}.

Since $\mu$ is invariant nonatomic, we obtain
\begin{eqnarray*}
\mu(\Phi_{\frac{\epsilon}2}(x))&=&\mu(\Phi_{\frac{\epsilon}2}(x)\setminus\{x\})\\
&\overset{\eqref{pajuo}}{\leq}& \sum_{i\geq0}\mu(T^{-i}(B(T^i(x),\epsilon)\setminus C(T^i(x))))\\
&=& \sum_{i\geq0}\mu(B(T^i(x),\epsilon)\setminus C(T^i(x)))\\
&\overset{\eqref{pa}}{=}&0
\end{eqnarray*}
for all $x\in X$. Therefore, $\frac{\epsilon}2 $ is an expansivity constant of $\mu$ and so $\mu$ is positively expansive.
This completes the proof.
\end{proof}

\begin{proof}[Proof of Theorem \ref{ermantraut}]
Let $T:X\to X$ be a continuous map of a compact metric space.
We shall prove that every ergodic invariant measure $\mu$ of $T$
satisfies \eqref{jimmy}.

If $E_\mu(T)\leq 0$, then the right-hand side of \eqref{jimmy} is nonpositive and we are done since the Kolmogorov-Sinai entropy is nonnegative.

Now, assume that $E_\mu(T)>0$.
Fix $0<\Delta<E_\mu(T)$.
Then, there are $\lambda>E_\mu(T)-\Delta$ and
$\epsilon>0$ such that \eqref{pito} holds this $\epsilon$ and
$k=e^\lambda$.

Given $x\in X$, $n\in\mathbb{N}$ and $\gamma>0$ we define
$$
B(x,n,\gamma)=\{y\in X\mid d(T^i(x),T^i(y))<\gamma,\,\forall i=0,\cdots, n-1\}
$$
and if $F\subset X$,
$$
B(F,n,\gamma)=\bigcup_{x\in F}B(x,n,\gamma).
$$

For all $x\in X$ we define $C(x)$ as in \eqref{tos} hence \eqref{pa} holds.
Since $\mu$ is invariant,
\begin{equation}
\label{eso}
\mu\left(
\bigcup_{i=0}^{n-1} T^{-i}(B(T^i(x),\gamma)\setminus C(T^i(x)))\right)=0,
\quad\quad\forall x\in X,\,\forall n\in\mathbb{N},\,0<\gamma<\epsilon.
\end{equation}
Now, fix $x\in X$.
Take an integer $n\geq2$, $0<\gamma<\epsilon$ and
$$
y\in B(x,n,\gamma)\setminus \left(
\bigcup_{i=0}^{n-1} T^{-i}(B(T^i(x),\gamma)\setminus C(T^i(x)))\right).
$$
It follows that
$$
T^i(y)\in B(T^i(x),\gamma)\quad(\forall0\leq i\leq n-1)\quad\mbox{ and }\quad
y\notin \bigcup_{i=0}^{n-1} T^{-i}(B(T^i(x),\gamma)\setminus C(T^i(x))).
$$
From both properties we obtain
$$
y\in \bigcap_{i=0}^{n-1} T^{-i}(C(T^i(x))).
$$
In particular,
$T^{n-2}(y)\in C(T^{n-2}(x))$ so
$$
d(T^{n-1}(x),T^{n-1}(y))\geq kd(T^{n-2}(x),T^{n-2}(y)).
$$
Likewise,
$T^{n-3}(y)\in C(T^{n-3}(x))$ so
$$
d(T^{n-2}(x),T^{n-2}(y))\geq kd(T^{n-3}(x),T^{n-3}(y)).
$$
Repeating this process we obtain
$$
d(T^{n-i}(x),T^{n-i}(y))\geq kd(T^{n-i-1}(x),T^{n-i-1}(y)),\quad\quad\forall
1\leq i\leq n-1.
$$
Since $d(T^{n-1}(x),T^{n-1}(y))<\gamma$, we obtain
$\gamma>k^{n-1}d(x,y)$ so $y\in B(x,k^{-(n-1)}\gamma)$.
We conclude that
\begin{equation}
\label{cuentas}
B(x,n,\gamma)\setminus \left(
\bigcup_{i=0}^{n-1} T^{-i}(B(T^i(x),\gamma)\setminus C(T^i(x)))\right)\subset B(x,k^{-(n-1)}\gamma),
\end{equation}
for all $x\in X$, $n\in\mathbb{N}$ and $0<\gamma<\epsilon.$

On the other hand, $\mu$ is ergodic invariant and $T$ is continuous so Katok $\delta$-entropy formula (Theorem 1.1 in \cite{k}) implies
\begin{equation}
\label{katok}
h_\mu(T)=\lim_{\gamma\to0}\limsup_{n\to\infty}\frac{1}n\log r_\mu(n,\gamma,\delta), \quad\quad\forall \delta>0,
\end{equation}
where $r_\mu(n,\gamma,\delta)$ is the minimal cardinality $card(F)$ of a {\em $\mu$-$(n,\gamma,\delta)$-spanning set} i.e. $F\subset X$ such that
$\mu(B(F,n,\gamma))\geq 1-\delta$.

Let $F$ be a finite $\mu$-$(n,\gamma,\delta)$-spanning set for some
$n\in\mathbb{N}$, $\delta>0$ and $0<\gamma<\epsilon$
(it exists since $X$ is compact and $B(x,n,\gamma)$ is an open neighborhood of $x$ for all $n\in\mathbb{N}$ and $\epsilon>0$).

It follows that the collection $\{B(x,n,\gamma)\mid x\in F\}$ covers a set of $\mu$-measure $\geq 1-\delta$.
Since $F$ is finite, 
\eqref{eso} implies
$$
\mu\left(
\bigcup_{x\in F}\bigcup_{i=0}^{n-1} T^{-i}(B(T^i(x),\gamma)\setminus C(T^i(x))
\right)=0.
$$
Then,
\begin{eqnarray*}
1-\delta&\leq& \mu\left(\bigcup_{x\in F}B(x,n,\gamma)\right)\\
&=& \mu\left(
\bigcup_{x\in F}B(x,n,\gamma)\setminus \bigcup_{x\in F}\bigcup_{i=0}^{n-1} T^{-i}(B(T^i(x),\gamma)\setminus C(T^i(x)))
\right)\\
&\leq&
\mu\left(
\bigcup_{x\in F}\left[B(x,n,\gamma)\setminus\bigcup_{i=0}^{n-1} T^{-i}(B(T^i(x),\gamma)\setminus C(T^i(x)))\right]
\right)\\
&\overset{\eqref{cuentas}}{\leq}& \mu\left(\bigcup_{x\in F}B(x,k^{-(n-1)}\gamma)\right)
\end{eqnarray*}
and so the collection
of balls $\{B(x,k^{-(n-1)}\gamma)\mid x\in F\}$ also covers a set of $\mu$-measure $\geq 1-\delta$.

This implies $N_\mu(X,k^{-(n-1)}\gamma,\delta)\leq card(F)$ and, since $F$ is an arbitrary finite $\mu$-$(n,\gamma,\delta)$-spanning set, we obtain
$$
N_\mu(X,k^{-(n-1)}\gamma,\delta)\leq r_\mu(n,\gamma,\delta),\quad\quad\forall n\geq N,\,\delta>0,\, 0<\gamma<\epsilon.
$$
Then,
$$
\limsup_{n\to\infty}\frac{\log N_\mu(X,k^{-(n-1)}\gamma,\delta)}{-\log(k^{-(n-1)}\gamma)}\leq \frac{1}{\log k}\limsup_{n\to\infty}\log r_\mu(n,\gamma,\delta),\quad\forall \delta>0,\, 0<\gamma<\epsilon.
$$
But the map $\varphi:(0,\infty)\to\mathbb{R}$ defined by
$\varphi(\beta)=N_\mu(X,\beta,\delta)$ is decreasing for all $\delta>0$ so
$$
\limsup_{\beta\to0}\frac{N_\mu(X,\beta,\delta)}{-\log\beta}=
\limsup_{n\to\infty}\frac{\log N_\mu(X,k^{-(n-1)}\gamma,\delta)}{-\log(k^{-(n-1)}\gamma)}
$$
by Lemma 2.6 in Fathi \cite{f1} thus
$$
\limsup_{\beta\to0}\frac{N_\mu(X,\beta,\delta)}{-\log\beta}\leq\frac{1}{\log k}\limsup_{n\to\infty}\log r_\mu(n,\gamma,\delta),\quad\quad\forall\delta>0,\,0<\gamma<\epsilon.
$$
Letting $\gamma\to0$ above using \eqref{katok} we get
$$
\limsup_{\beta\to0}\frac{N_\mu(X,\beta,\delta)}{-\log\beta}\leq \frac{h_\mu(T)}{\log k},\quad\quad\forall \delta>0.
$$
Letting $\delta\to0$ above using \eqref{faraday} we get
$$
\overline{dim}_B(\mu)\leq \frac{h_\mu(T)}{\log k}.
$$
But $k=e^\lambda$ with $\lambda>E_\mu(T)-\Delta>0$ so
$$
\overline{dim}_B(\mu)\cdot (E_\mu(T)-\Delta)<h_\mu(T),\quad\quad\forall0<\Delta<E_\mu(T).
$$
Letting $\Delta\to0$ above we get \eqref{jimmy} completing the proof.
\end{proof}

\section*{DElaration of competing interest}

\noindent
There is no competing interest.

\section*{Data availability}

\noindent
No data was used for the research described in the article.

\end{document}